\documentclass[twoside,twocolumn]{article}

\usepackage{blindtext} 

\usepackage[english]{babel} 

\usepackage{titling} 

\usepackage{hyperref} 

\usepackage{amsmath}

\title{The minimum overlap problem revisited} 
\author{%
\textsc{Jan Kristian Haugland}\\
\normalsize \href{mailto:admin@neutreeko.net}{admin@neutreeko.net} 
}
\date{\today} 

\renewcommand{\maketitlehookd}{%
	\begin{abstract}
		\noindent For a given partition of (1, 2, ..., 2n) into two disjoint subsets A and B with n elements in each, consider the maximum number of times any integer occurs as the difference between an element of A and an element of B. The minimum value of this maximum (over all partitions) is denoted by M(n). By a result of Swinnerton-Dyer, one way to estimate lim M(n)/n from above is to give step functions that describe the density of A, say, throughout the interval [1, 2n] for a large n rather than looking for explicit partitions. A step function that improves the upper bound from 0.382002... to 0.380926... is given.
	\end{abstract}
}

\begin{document}

\maketitle


Consider a partition of $\lbrace 1, 2, ..., 2n\rbrace $ into two disjoint subsets $\lbrace a_i \rbrace $ and $\lbrace b_j \rbrace $ with $n$ elements in each. For a fixed integer $k$, denote by $M_k$ the number of solutions to $a_i - b_j = k$, and let $M(n)$ denote the minimum, over all partitions, of $\max_k M_k$. To estimate $M(n)$ is the minimum overlap problem of Paul Erd\"{o}s.

Swinnerton-Dyer proved in \cite{Haugland:1996} that $\lim_{n\to \infty} {M(n)\over n}$ is equal to the infimum, over all step functions $f$ on $\left[ 0, 2\right]$ with values in $\left[ 0, 1\right]$ and satisfying $$\int_0^2 f(x)dx = 1$$ of \begin{equation} \label{eq:1} \max_k \int f(x)(1-f(x+k))dx \end{equation}
For simplicity, we let "a step function with $n$ steps" denote a function that is constant on any interval $\left( {2i\over n}, {2(i+1)\over n}\right)$ where $i \in \lbrace 0, 1, ..., n - 1 \rbrace$. A step function with 21 steps for which (1) attains the value 0.382002... is given in the same paper. The purpose of this note is to present an improvement on this example. In comparison, the best known lower bound for $\lim_{n\to \infty} {M(n)\over n}$ is $\sqrt{4-\sqrt{15}}=0.356393...$ by Leo Moser \cite{Moser:1959}.

An improvement on the upper bound can be found using a step function with only 15 steps. Taking
$$f(x)=0 \text{ for } 0\leq x < {2\over 15}$$
$$f(x)=0.09938602 \text{ for } {2\over 15}\leq x < {4\over 15}$$
$$f(x)=0.64299877 \text{ for } {4\over 15}\leq x < {6\over 15}$$
$$f(x)=0.36104582 \text{ for } {6\over 15}\leq x < {8\over 15}$$
$$f(x)=0.69536426 \text{ for } {8\over 15}\leq x < {10\over 15}$$
$$f(x)=0.59241335 \text{ for } {10\over 15}\leq x < {12\over 15}$$
$$f(x)=0.89573331 \text{ for } {12\over 15}\leq x < {14\over 15}$$
$$f(x)=0.92611694 \text{ for } {14\over 15}\leq x \leq 1$$
$$f(x)=f(2-x) \text{ for } 1 < x \leq 2$$
yields the value 0.38153155 (when rounded upwards) for (1). Using 19 steps allows for a further improvement to the value 0.381112263316104816. This is attained by taking
$$f(x)=0 \text{ for } 0\leq x < {4\over 19}$$
$$f(x)=0.348795091509472207 \text{ for } {4\over 19}\leq x < {6\over 19}$$
$$f(x)=0.742684181900847446 \text{ for } {6\over 19}\leq x < {8\over 19}$$
$$f(x)=0.207655267155520404 \text{ for } {8\over 19}\leq x < {10\over 19}$$
$$f(x)=0.780222086674911898 \text{ for } {10\over 19}\leq x < {12\over 19}$$
$$f(x)=0.568104573396874436 \text{ for } {12\over 19}\leq x < {14\over 19}$$
$$f(x)=0.689049157609512654 \text{ for } {14\over 19}\leq x < {16\over 19}$$
$$f(x)=0.967251286500411737 \text{ for } {16\over 19}\leq x < {18\over 19}$$
$$f(x)=0.892476710504898436 \text{ for } {18\over 19}\leq x \leq 1$$
$$f(x)=f(2-x) \text{ for } 1 < x \leq 2$$
The best upper bound we have found comes from a step function with 51 steps. Taking
$$f(x)=0 \text{ for } 0\leq x < {10\over 51}$$
$$f(x)=0.0002938681556273 \text{ for } {10\over 51}\leq x < {12\over 51}$$
$$f(x)=0.5952882223921177 \text{ for } {12\over 51}\leq x < {14\over 51}$$
$$f(x)=0.7844530825484313 \text{ for } {14\over 51}\leq x < {16\over 51}$$
$$f(x)=0.8950034338013842 \text{ for } {16\over 51}\leq x < {18\over 51}$$
$$f(x)=0.0597964076006748 \text{ for } {18\over 51}\leq x < {20\over 51}$$
$$f(x)=0.0189602838469592 \text{ for } {20\over 51}\leq x < {22\over 51}$$
$$f(x)=0.7420501628172980 \text{ for } {22\over 51}\leq x < {24\over 51}$$
$$f(x)=0.6444559588500921 \text{ for } {24\over 51}\leq x < {26\over 51}$$
$$f(x)=0.3549040817844764 \text{ for } {26\over 51}\leq x < {28\over 51}$$
$$f(x)=0.8762442385073478 \text{ for } {28\over 51}\leq x < {30\over 51}$$
$$f(x)=0.5437907313675501 \text{ for } {30\over 51}\leq x < {32\over 51}$$
$$f(x)=0.2679640048997296 \text{ for } {32\over 51}\leq x < {34\over 51}$$
$$f(x)=0.8518954615823791 \text{ for } {34\over 51}\leq x < {36\over 51}$$
$$f(x)=0.5211171156914872 \text{ for } {36\over 51}\leq x < {38\over 51}$$
$$f(x)=1 \text{ for } {38\over 51}\leq x < {40\over 51}$$
$$f(x)=0.5506146790047043 \text{ for } {40\over 51}\leq x < {42\over 51}$$
$$f(x)=0.9007715390796991 \text{ for } {42\over 51}\leq x < {44\over 51}$$
$$f(x)=0.8229000691941086 \text{ for } {44\over 51}\leq x < {46\over 51}$$
$$f(x)=0.8879541710440111 \text{ for } {46\over 51}\leq x < {48\over 51}$$
$$f(x)=0.9315424878319221 \text{ for } {48\over 51}\leq x < {50\over 51}$$
$$f(x)=1 \text{ for } {50\over 51}\leq x \leq 1$$
$$f(x)=f(2-x) \text{ for } 1 < x \leq 2$$
yields the value 0.3809268534330870 for (1).



\end{document}